\documentclass[12pt]{article}
\usepackage{amssymb}
\usepackage{amsfonts}
\usepackage{amsmath}
\usepackage[usenames]{color}
\usepackage{mathrsfs}
\usepackage{amsfonts}
\usepackage{amssymb,amsmath}
\usepackage{CJK}
\usepackage{cite}
\usepackage{cases}
\usepackage{amsthm}\usepackage{multirow}
\usepackage{graphicx}
\usepackage{float}
\usepackage{makecell}
\usepackage{latexsym}
\usepackage{CJK}
\usepackage{caption}
\usepackage{bm}
\usepackage [latin1]{inputenc}

\def\cl{\centerline}

\def\vs{\vspace*}

\def\ni{\noindent}

\numberwithin{equation}{section}
\newtheorem{theo}{Theorem}[section]

\newtheorem{lemm}[theo]{Lemma}
\newtheorem{exam}[theo]{Example}
\newtheorem{prop}[theo]{Proposition}

\newtheorem{case}{Case}

\newtheorem{rema}[theo]{Remark}

\newtheorem*{thmA}{Theorem A}
\newtheorem*{thmB}{Theorem B}

\newcommand{\arccot}{\mathrm{arccot}\,}
\begin{document}
\captionsetup[figure]{labelfont={bf},labelformat={default},labelsep=period,name={Fig.}}

\begin{center}
\cl{\large\bf \vs{6pt} The isoparametric functions on a class of Finsler spheres\,{$^*\,$}}
\footnote {$^*\,$ Project supported by AHNSF (No.2108085MA11).
\\\indent\ \ $^\dag\,$ hequn@tongji.edu.cn
}
\cl{Yali Chen$^1$, Qun He$^1$$^\dag\,$}

\cl{\small 1 School of Mathematical Sciences, Tongji University, Shanghai,
200092, China.}
\end{center}

{\small
\parskip .005 truein
\baselineskip 3pt \lineskip 3pt

\noindent{{\bf Abstract:}
In this paper, we give global expressions of geodesics and isoparametric functions on a Randers sphere by navigation. We obtain isoparametric families and focal submanifolds in $(\mathbb{S}^{n}, F, d\mu_{BH})$ by Cartan-M\"{u}nzner polynomials. Further more, we construct some examples of closed and non-closed geodesics, isoparametric functions, isoparametric families and focal submanifolds.

\ni{\bf Key words:}
 geodesic; isoparametric hypersurface; isoparametric function; Randers sphere.}

\ni{\it Mathematics Subject Classification (2010):} 53C60, 53C42, 34D23.}
\parskip .001 truein\baselineskip 6pt \lineskip 6pt
\section{Introduction}
~~~~The study on isoparametric hypersurfaces is an important topic in Riemannian geometry. E. Cartan began to study the isoparametric hypersurfaces in real space forms systematically\cite{C}. Then many mathematicians started working on it and made many important contributions~\cite{C1,TP,GT}. In Finsler geometry, the conception of isoparametric hypersurfaces has been introduced in~\cite{HYS}. Let $(N,F,d\mu)$ be an
$n$-dimensional Finsler manifold with volume form $d\mu$. A function $f$
on $(N,F,d\mu)$ is called \textit{isoparametric} if it is almost everywhere smooth and there are functions $\widetilde{a}(t)$ and $\widetilde{b}(t)$ such that
\begin{equation}\label{1.1} \left\{\begin{aligned}
&F(\nabla f)=\tilde{a}(f),\\
&\Delta f=\tilde{b}(f),
\end{aligned}\right.
\end{equation}
where $\nabla f$ denotes the gradient of $f$, which is defined by means of the Legendre transformation, and $\Delta f$ is a nonlinear Finsler-Laplacian of $f$ (See Section 2.2 for details). Each regular level hypersurface of an isoparametric function is called an isoparametric hypersurface. As in Riemannian geometry, We call a complete and simply connected Finsler manifold with constant flag curvature a Finsler space form. Minkowski spaces (with zero flag curvature) and Finsler spheres (with positive constant flag curvature) are two important classes of Finsler space forms~\cite{YH}. For a Finsler sphere $(\mathbb{S}^{n},F)$, if $F$ is a Randers metric, we call it a Randers sphere.~\cite{HYS,HYS1,HD,HDY} gave a complete classification of $d\mu_{BH}-$isoparametric hypersurfaces in Minkowski spaces, Randers space forms and Funk-type spaces.~\cite{HCY} studied isoparametric hypersurfaces in Finsler space forms by investigating focal points, tubes and parallel hypersurfaces.

The study on isoparametric hypersurfaces is inseparable from isoparametric functions. In Riemannian geometry, M\"{u}nzner showed that an isoparametric hypersurface in $\mathbb{S}^{n}\subseteq\mathbb{R}^{n+1}$ with $g$ distinct principal curvatures is obtained by a homogeneous polynomial of degree $g$ on $\mathbb{R}^{n+1}$ satisfying the Cartan-M\"{u}nzner differential equations\cite{M1,M}. In Finsler geometry, it is still an open problem to study isoparametric functions on $\mathbb{S}^{n}$. It has been proved that the isoparametric hypersurfaces in Euclidean spheres and Randers spheres are the same, but their isoparametric functions, isoparametric families and focal submanifolds are all different, except in very special cases \cite{HDY}. From \cite{ZHC}, we know the isoparametric family is a special mean curvature flow in Finsler manifolds, which means the mean curvature flows generated by the same isoparametric hypersurface are different in general. So it is natural and meaningful to study isoparametric functions on Randers spheres and their corresponding isoparametric families.

Geodesic is important in Finsler geometry, especially in the study of isoparametric theory. We have already known the geodesics in Minkowski spaces and Funk spaces are all straight lines. But geodesics on a Randers sphere are more complicated. They may not be a closed great circle. C. Robles classified geodesics on a Randers manifold of constant flag curvature~\cite{RO}. M. Xu proved that the number of geometrically distinct closed geodesics on $(\mathbb{S}^n, F)$ is at least dim $I(\mathbb{S}^n, F)$\cite{XM}.

Navigation is a technique to manufacture new Finsler metrics and study their geometric properties. D. Bao, C. Robles and Z. M. Shen classified Randers metrics of constant flag curvature by navigation on Riemannian manifolds\cite{BRS}. L. Huang and X. Mo gave a geometric description of the geodesics of the Finsler metric produced from any Finsler metric and any homothetic field in terms of a navigation presentation\cite{HM12,HM1}. M. Xu, V. S. Matveev, et al discussed the correspondences of geodesics and Jacobi fields for homothetic navigation\cite{XM1}.

The local correspondences of isoparametric functions for homothetic navigation have been given in \cite {XM1}. In this paper, we try to give a global expression of isoparametric functions on a Randers sphere using navigation skills via a global expression of geodesics. Further more, we try to obtain the corresponding isoparametric families and focal submanifolds.

Let $(\mathbb S^n, F_{Q})$ be a Randers sphere, where $F_{Q}$ is a Randers metric with the navigation datum $(h, V)$, $h$ is a standard Euclidean sphere metric and $V=Qx\mid_{\mathbb S^{n}}$,$~x\in \mathbb R^{n+1}$, $Q\in o(n+1)$ such that $|V|<1$. We can prove that $V=Qx\mid_{\mathbb S^{n}}$ is a Killing vector field with respect to $h$. In fact, any Killing vector field on $(\mathbb S^{n},h)$ can be obtained in this way. Hence, we have

\begin{theo} \label{thm01}
Let $\overline{f}$ be a global isoparametric function on $(\mathbb S^n, h)$. Set $|\nabla^{h}\overline{f}|=a(\overline{f})$, $\overline{f}(\mathbb{S}^n)=[c,d]$, then $\psi(x)=\exp(Q\int^{\overline{f}(x)}_{t_{0}}\frac{1}{a(t)}dt)x$ is a homeomorphism from $\mathbb{S}^n$ to $\mathbb{S}^n$ and $f=\overline{f}\circ\psi^{-1}$ is a global isoparametric function on $(\mathbb S^n, F_{Q})$. Conversely, for any isoparametric function $f$ on $(\mathbb S^n, F_{Q})$ which satisfies $F_{Q}(\nabla f)=a(f)$, where $a^{2}(t)\in C^{2}(f(\mathbb{S}^n))$, can be obtained in this way.
\end{theo}
In Riemannian case, M\"{u}nzner made the association between an isoparametric function and a homogeneous polynomial.
\begin{thmA}~\cite{TP}
Let $M\subset\mathbb S^{n}\subset\mathbb R^{n+1}$ be a connected isoparametric hypersurface with $g$ principal curvatures $\lambda_{i}=cot\theta_{i}$, $0<\theta_{i}<\pi$, with multiplicities $m_{i}$. Then $M$ is an open subset of a level set of the restriction to $\mathbb S^{n}$ of a homogeneous polynomial $\phi$ on
$\mathbb R^{n+1}$ of degree $g$ satisfying the differential equations,
\begin{equation}\label{1.4}
\left\{\begin{array}{ll}
|\nabla^E\phi|^2=g^{2}r^{2g-2},\\
\Delta^{E}\phi=cr^{g-2},
\end{array}\right.
\end{equation}
where $r=|x|$, and $c=\frac{g^{2}(m_{2}-m_{1})}{2}$, $m_{1}$, $m_{2}$ are the two distinct multiplicities. $\phi$ is called the Cartan-M\"{u}nzner polynomial of $M$. (\ref{1.4}) is called the Cartan-M\"{u}nzner differential equations.
\end{thmA}
\begin{thmB}~\cite{TP}
Let $\phi: \mathbb{R}^{n+1}\rightarrow\mathbb{R}$ be a Cartan-M\"{u}nzner polynomial of degree $g$ and $f$ is the restriction to $\mathbb{S}^{n}$. Then each isoparametric hypersurface
$$M_{t}=f^{-1}(t),\ \ \ \ -1<t<1$$
is connected. Moreover, $M_{+}=f^{-1}(1)$ and $M_{-}=f^{-1}(-1)$ are the focal submanifolds, respectively, and they are also connected.
\end{thmB}
Due to Theorem A, Theorem B and Theorem \ref{thm01}, we can characterize all isoparametric families of isoparametric hypersurfaces in a Randers sphere by a Cartan-M\"{u}nzner polynomial.
\begin{theo} \label{thm02}
Let $\{M_{t}\}$ be a connected isoparametric family in $(\mathbb S^{n}, F_{Q})$ with $g$ principal curvatures with multiplicities $m_{i}$. Then there exists a Cartan-M\"{u}nzner polynomial $\phi$ of degree $g$ such that $f=\phi\circ\psi^{-1}\mid_{\mathbb S^{n}}$ is an isoparametric function on $(\mathbb S^{n}, F_{Q})$, where $\psi(x)=\exp\left(\frac{1}{g}\arcsin\frac{\phi(x)}{|x|^{g}}Q\right)x$, and each $M_{t}$ is an open subset of a level set of $f$. Further more, the maximized connected isoparametric family and two connected focal submanifolds can be expressed as
$$f^{-1}(t)=\{\exp(\frac{\arcsin t}{g}Q)((\cos\frac{\arcsin t}{g})x+\sin(\frac{\arcsin t}{g})(\mathbf{n}-V))\mid x\in M\},\ \ -1<t<1$$
and
$$M_{\pm}=\{\exp(\pm \frac{\pi}{2g}Q)(\cos\frac{\pi}{2g}x\pm \sin\frac{\pi}{2g}(\mathbf{n}-V))\mid x\in M\},$$
where $M=f^{-1}(0)$ and $\mathbf{n}$ is the unit normal vector of $M$ at $x$.
\end{theo}
\begin{rema}
Every isoparametric family is a special mean curvature flow. Hence, from a given isoparametric hypersurface $M$, we can construct two different mean curvature flows in $(\mathbb{S}^{n},h)$ and $(\mathbb{S}^{n},F_{Q})$ by a Cartan-M\"{u}nzner polynomial, respectively.
\end{rema}
The contents of this paper are organized as follows. In Section 2, some fundamental concepts and formulas are given. In Section 3,
we give global expressions of a one-parameter group and geodesics on $(\mathbb S^{n},F_{Q})$. Particularly, when $n=2$, the images of all geodesics can be depicted concretely. In section 4, we characterize isoparametric functions on $(\mathbb S^{n},F_{Q})$ and obtain global expressions of isoparametric families and focal submanifolds by homogeneous polynomial functions in  $\mathbb{R}^{n+1}$. Further more, we give some examples of isoparametric functions, isoparmetric families and focal submanifolds.
\section{Preliminaries}
~~~~In this section, we will give some definitions and lemmas that will be used in the proof of our main results.
\subsection{Finsler manifolds}
~~~~Let $N$ be a manifold and let $TN=\cup_{x\in N}T_xN$ be the tangent bundle of $N$, where $T_xN$ is the tangent space at
$x\in N$. A Finsler metric is a Riemannian metric without quadratic restriction. Precisely, a function $F(x,y)$ on $TN$ is called a Finsler metric on a manifold $N$ with local coordinates $(x,y)$, where $x=(x^i)$ and $y=y^i\frac{\partial}{\partial x^{i}}$ ,
if it has the following properties:

(i)\ \ Regularity:\ \ $F(x,y)$ is $C^{\infty}$ on $TN\backslash\{0\}$;

(ii)\ \ Positive homogeneity:\ \ $F(ty)=tF(y),\ \forall t>0, y\in T_xN$;

(iii)\ \ Strong convexity:\ \ The $n\times n$ matrix $(\frac{\partial^2F^2}{\partial y^i \partial y^j}(x,y))(y\neq 0)$ is positive definite.

The fundamental form~$g$ of~$(N,F)$ is
\begin{equation*}
g=g_{ij}(x,y)dx^{i} \otimes dx^{j}, ~~~~~~~g_{ij}(x,y)=\frac{1}{2}[F^{2}] _{y^{i}y^{j}}.
\end{equation*}
For a Finsler metric $F=F(x,y)$ on a manifold $N$, the geodesic $\gamma=\gamma(t)$ of $F$ in local coordinates $(x^i)$ are characterized by
$$\frac{\textmd{d}^2x^i}{\textmd{d}t^2}+2G^i(x,\frac{\textmd{d}x}{\textmd{d}t})=0,$$
where $(x^i(t))$ are the coordinates of $c(t)$ and $G^i=G^i(x,y)$ are defined by
$$G^i=\frac{g^{il}}{4}\{[F^2]_{x^ky^l}y^k-[F^2]_{x^l}\},$$
which are called \textit{the spray coefficients}. For $X=X^{i}\frac{\partial}{\partial x^{i}}\in\Gamma(TN)$, the covariant derivative of $X$ along $v=v^{i}\frac{\partial}{\partial x^{i}}\in T_{x}N$ with respect to a reference vector $w\in T_{x}N\setminus0$ is defined by
$$\nabla^{w}_{v}X(x):=\{v^{j}\frac{\partial X^{i}}{\partial x^{j}}(x)+\Gamma^{i}_{jk}(w)v^{j}X^{k}(x)\}\frac{\partial}{\partial x^{i}}.$$
The equation of geodesics can be expressed by  $\nabla^{\dot{\gamma}}_{\dot{\gamma}}\dot{\gamma}\equiv0$.

Let~${\mathcal L}:TN\rightarrow T^{\ast}N$ denote the \emph{Legendre transformation}, satisfying~${\mathcal L}(\lambda
y)=\lambda {\mathcal L}(y)$ for all~$\lambda>0,~y\in TN$.
For a smooth function~$f: N\rightarrow \mathbb{R}$, the \emph{gradient vector} of~$f$ at~$x$ is defined as~$\nabla f(x)={\mathcal
L}^{-1}(df(x))\in T_{x}N$.
Set~$N_{f}=\{x\in N|df(x)\neq 0\}$ and~$\nabla^{2}f(x)=\nabla^{\nabla f}(\nabla f)(x)$ for~$x\in N_{f}$. The \emph{Laplacian} of~$f$ can defined
by
\begin{equation}\label{2.1}
\hat{\Delta} f=\textmd{tr}_{g_{_{\nabla f}}}(\nabla^{2}f).
\end{equation}
And the \emph{Laplacian} of~$f$ with respect to the volume form~$d\mu=\sigma(x)dx=\sigma(x)dx^{1}\wedge dx^{2}\wedge\cdots\wedge dx^{n}$ can be represented as
\begin{equation}\label{2.2}
\Delta_{\sigma} f=\textmd{div}_{\sigma}(\nabla f)=\frac{1}{\sigma}\frac{\partial}{\partial x^{i}}(\sigma g^{ij}(\nabla f)f_{j})=\hat{\Delta} f-\textbf{S}(\nabla f),
\end{equation}
where
$$\textbf{S}(x,y)=\frac{\partial G^{i}}{\partial y^{i}}-y^{i}\frac{\partial}{\partial x^{i}}(\ln \sigma(x))$$
is the \emph{$\mathbf{S}$-curvature}~\cite{SZ}.
\subsection{Isoparametric functions and isoparametric hypersurfaces}
~~~~Let $f$ be a non-constant $C^1$ function defined on a Finsler manifold $(N, F, d\mu)$ and smooth on $N_{f}$. Set $J=f(N_{f})$. The function $f$ is called \textit{isoparametric  (resp. $d\mu-$isoparametric)} on $(N, F, d\mu)$, where $d\mu=\sigma(x)dx$, if there exist a smooth function $a(t)$ and a continuous function $b(t)$ defined on $J$ such that (\ref{1.1}) holds for $\Delta f=\hat{\Delta} f$ (resp. $\Delta f=\Delta_{\sigma}f$), which is defined by (\ref{2.1}) (resp.(\ref{2.2})). All the regular level hypersurfaces $M_{t}=f^{-1}(t)$ are named an \textit{($d\mu$-)isoparametric family}, each of which is called an \textit{($d\mu$-)isoparametric hypersurface} in $(N, F, d\mu)$. If each $M_{t}$ is connected, it is a connected isoparametric family.

Since $(\mathbb{S}^n,F_{Q}, d\mu_{BH})$ is a Randers sphere with $\textbf{S}=0$, then $M$ is an isoparametric hypersurface if and only if $M$ is a $d\mu_{BH}$-isoparametric hypersurface.
\subsection{Geometric correspondences for homothetic navigation}
~~~~Let $V$ be a smooth vector field on a Finsler manifold $(N, F)$. Around each $x\in N$, $V$ generates a family of local diffeomorphisms $\psi_{t}$, which is a flow in an alternative terminology and satisfies
\begin{equation}\label{0.3}
\left\{\begin{array}{ll}
\psi_{0}=id: N\rightarrow N,\\
\psi_{s}\circ\psi_{t}=\psi_{s+t}, \ \ \ \ \forall s, t\in(-\varepsilon, \varepsilon), s+t\in(-\varepsilon, \varepsilon).
\end{array}\right.
\end{equation}
$V$ is called a \textit{homothetic field} when it satisfies
$$(\psi_{t}^{\ast}F)(x,y)=F(\psi_{t}(x),\psi_{t*}(y))=e^{-2ct}F(x,y),$$
where $x\in N$, $y\in T_{x}N$ and $t\in\mathbb{R}$. The constant $c$ is the dilation of $V$. If $c=0$, $V$ is Killing~\cite{SS}.
The main technique of the navigation problem is described as follows. Suppose $F$ is a Finsler metric and $V$ is a vector field with $F(x,-V)<1$, we can define a new Finsler metric $\tilde{F}$ by
$$F(x, y-\tilde{F}(x,y)V)=\tilde{F}(x,y), \ \ \ \ \ \ \forall x\in N,\ \ y\in T_{x}N.$$
\begin{lemm} \label{lemm33}~\cite{HM12}
Let $F=F(x,y)$ be a Finsler metric on a manifold $N$ and let $V$ be a vector field on $N$ with $F(x, -V_{x})<1$. Suppose that $V$ is homothetic with dilation $c$. Let $\tilde{F}=\tilde{F}(x,y)$ denote the Finsler metric on $N$. Then the geodesics of $\tilde{F}$ are given by $\psi_{t}(\gamma(a(t)))$, where $\psi_{t}$ is the flow of $V$, $\gamma(t)$ is a geodesic of $F$ and $a(t)$ is defined by
$$a(t):=\left\{\begin{array}{ll}
\frac{e^{2ct}-1}{2c},\ \ \ \ if\ \ c\neq0,\\
t,\ \ \ \ \ \ \ \ \ \ if\ \  c=0.
\end{array}\right.$$
\end{lemm}
If $F$ is a Randers metric with the navigation datum $(h, V)$, where $h$ is a Riemannian metric and $V$ is a Killing vector field, then we have
\begin{lemm}\label{lemm334}
For any geodesic $\overline{\gamma}(t)$ with respect to the metric $h$ satisfying $\overline{\gamma}(0)=\overline{x}$, $\gamma(t+t_{0})=\psi_{t+t_{0}}(\overline{\gamma}(t))$ is a geodesic with respect to the metric $F$ satisfying $\gamma(t_{0})=x=\psi_{t_{0}}(\overline{\gamma}(0))$.
\end{lemm}
\proof By ($\ref{0.3}$), $\gamma(t+t_{0})=\psi_{t+t_{0}}(\overline{\gamma}(t))=\psi_{t}(\psi_{t_{0}}(\overline{\gamma}(t)))$. Since $V$ is Killing, then $\psi_{t}$ is isometric with respect to the metric $h$. Hence, $\psi_{t_{0}}(\overline{\gamma}(t))$ is still a geodesic for the metric $h$ with $x=\psi_{t_{0}}(\overline{\gamma}(0))=\gamma(t_{0})$. Combine Lemma \ref{lemm33}, we complete the proof.
\endproof
\begin{lemm} \label{lemm46}~\cite{XM1}
Let $\widetilde{F}$ be the Finsler metric defined by navigation from the datum $(F,V)$ in which $V$ is a Killing vector field. Assume $x_{0}$ is a point where $F(x_{0},-V(x_{0}))<1$. Then for any normalized isoparametric function $f$ for $(F,d\mu^{F}_{BH})$ around the point $x_{0}$, the function $\widetilde{f}$ defined by $\widetilde{f}^{-1}(t)=\psi_{t}(f^{-1}(t))$ is a normalized isoparametric function for $(\widetilde{F},d\mu_{BH}^{\widetilde{F}})$ around $x_{0}$.
\end{lemm}
\section{ The geodesics of a Randers sphere}
\begin{lemm} \label{lemm31}
Let $x\in\mathbb{R}^{n+1}$, $Q\in o(n+1)$ and $h$ be a standard Euclidean sphere metric, then $V=Qx\mid_{\mathbb S^{n}}$ is a global Killing vector field on $(\mathbb S^n, h)$, and $|V|<1$ if and only if $I+Q^2$ is positive definite.
\end{lemm}
\proof Set $\widetilde{V}=Qx$, take $X,Y\in T\mathbb{R}^{n+1}$, we have
$$\langle\nabla_{X}^{E}\widetilde{V},Y\rangle=\langle Q\nabla_{X}^{E}x,Y\rangle=\langle QX,Y\rangle=-\langle X,QY\rangle=\langle X, Q\nabla_{Y}^{E}x\rangle=-\langle X, \nabla_{Y}^{E}\widetilde{V}\rangle,$$
where $\nabla^{E}$ is the covariant derivative in a Euclidean space. Hence, $\widetilde{V}$ is a global Killing vector field on $\mathbb{R}^{n+1}$.

Due to $\langle Qx, x\rangle=0$, we know
$V=\widetilde{V}\mid_{\mathbb S^{n}}=Qx\mid_{\mathbb S^{n}}\in T\mathbb{S}^n$. Set $X,Y\in T\mathbb{S}^n$, we have
$$h(\nabla^{h}_{X}V, Y)=\langle\nabla^{E}_{X}\widetilde{V}, Y\rangle=-\langle\nabla^{E}_{Y}\widetilde{V}, X\rangle=-h(\nabla^{h}_{Y}V, X).$$
Hence, $V=Qx\mid_{\mathbb S^{n}}$ is a global Killing vector field on $(\mathbb S^n, h)$.

Further more, $|V|<1$ if and only if $x^TQ^TQx<1=x^Tx$, which means $x^T(I+Q^2)x>0$. This completes the proof.
\endproof

\begin{lemm} \label{lemm32}
Suppose that $\psi_{t}(x_{0})$ is the integral curve of $V=Qx\mid_{\mathbb{S}^n}$ through $x(0)=x_{0}\in\mathbb{S}^{n}$, where $t\in\mathbb{R}$. Then $\psi_{t}=\exp(tQ)\mid_{\mathbb S^n}$.
\end{lemm}
\proof  Let $\widetilde{\psi}_{t}(x_{0})$ be the integral curve of $\widetilde{V}=Qx$, $\widetilde{x}(0)=x_{0}$. Due to the definition of the one-parameter local group,
$$\widetilde{x}'(t)=\widetilde{V}(\widetilde{x}(t))=Q\widetilde{x}(t).$$
Then
$$\widetilde{x}(t)=\exp (tQ)C,$$
where $C$ is a constant vector. Noting that $\widetilde{x}(0)=x_{0}$, we have
$$\widetilde{x}(t)=\widetilde{\psi}_{t}(x_{0})=\exp (tQ)x_{0}.$$
Further more,
$$|\widetilde{x}(t)|^2=x_{0}^{T}\exp (tQ)^{T}\exp (tQ)x_{0}=|x_{0}|^2=1,$$
we know that $\widetilde{x}(t)$ is still on $\mathbb{S}^{n}$. From the existence and uniqueness theorem of initial value problem of ODE, $\widetilde{\psi}_{t}(x_{0})=\widetilde{x}(t)=\psi_{t}(x_{0})\in\mathbb{S}^{n}$ is the unique integral curve through $x_{0}$. Since $\widetilde{\psi}_{t}$ is a global one-parameter group on $\mathbb{R}^{n+1}$, then $\psi_{t}$ is a global one-parameter group on $\mathbb{S}^{n}$.
\endproof
 Let $(\mathbb S^n, F_{Q})$ be a Randers sphere, where $F_{Q}$ is a Randers metric with navigation datum $(h, V)$, $h$ is a standard Euclidean sphere metric, $V=Qx\mid_{\mathbb S^{n}},~x\in \mathbb R^{n+1}$, $Q\in o(n+1)$ such that $I+Q^2>0$. Then we have
\begin{theo} \label{thm0}
The unit speed geodesic $\gamma(s)$ on $(\mathbb S^n, F_{Q})$ satisfying $\gamma(0)=x\in\mathbb S^n$ and $\dot{\gamma}(0)=X\in T_x\mathbb S^n$ can be expressed as
$$\gamma(s)=\exp (sQ)((\cos s)x+(\sin s)(X-V(x))).$$
\end{theo}
\proof From Lemma \ref{lemm31} and Lemma \ref{lemm32}, we know that there exist a global Killing vector field $V=Qx\mid_{\mathbb{S}^n}$ and flow of diffeomorphisms $\psi_{t}=\exp (tQ)\mid_{\mathbb{S}^n}$ on $\mathbb{S}^{n}$. Let $\overline{\gamma}(s)$ be a unit speed geodesic on $(\mathbb S^n, h)$ satisfying $\overline{\gamma}(0)=x$ and $\dot{\overline{\gamma}}(0)=\overline{X}$, then $$\overline{\gamma}(s)=(\cos s)x+(\sin s)\overline{X},$$ where $s\in\mathbb{R}$.  From Lemma \ref{lemm33},
$$\gamma(s)=\psi_{s}(\overline{\gamma}(s))=\exp (sQ)(\overline{\gamma}(s))$$
is a local unit speed geodesic on $(\mathbb{S}^n,F_{Q})$ satisfies $\gamma(0)=x$, $\dot{\gamma}(0)=X$, where $s\in (-\varepsilon, \varepsilon)$. Obviously, $\psi_{t}$ satisfies (\ref{0.3}).

From Lemma \ref{lemm334}, we know that no matter what is the origin point of $\overline{\gamma}(t)$, $\gamma(t)$ is a geodesic on $(\mathbb S^n, F_{Q})$. Hence, $\gamma(s)=\psi_{s}(\overline{\gamma}(s))$ is a global geodesic on $(\mathbb{S}^n,F_{Q})$, where $s\in \mathbb{R}$. By the notion of navigation, $X=\dot{\gamma}(0)=\overline{X}+V$. Hence, $$\gamma(s)=\exp (sQ)(\overline{\gamma}(s))=\exp (sQ)((\cos s)x+\sin s(X-V)),~~s\in\mathbb{R}.$$  This completes the proof of Theorem \ref{thm0}.
\endproof

By Theorem \ref{thm0}, we can obtain all geodesics on $(\mathbb{S}^{n},F_{Q})$ and find out if the geodesics are closed or non-closed.
\begin{exam}\label{exam3}
In $(\mathbb{S}^2, F_{Q})$, we choose $x=(1,0,0)^{T}$, $\overline{X}=(0,1,0)^{T}$, then $\overline{\gamma}(s)=(\cos s, \sin s, 0)^{T}$. Take
\begin{equation}\label{3.1.1.1}
Q=\left(
\begin{array}{cccc}
0&a&b\\
-a&0&c\\
-b&-c&0
\end{array}
\right),
\end{equation}
where $I+Q^2$ is positive definite. By a direct calculation,
\begin{equation}\label{3.111}
\gamma(s)=\left(
\begin{array}{cc}
\cos bs\cos(1-a)s\\
\cos cs\sin(1-a)s-\sin bs\sin cs\cos(1-a)s\\
-\sin bs\cos cs\cos(1-a)s-\sin(1-a)s\sin cs
\end{array}
\right).
\end{equation}
Due to (\ref{3.111}), we try to find some closed or non-closed geodesics on $(\mathbb{S}^{2},F_{Q})$.
\begin{case}
The geodesics are closed if and only if there exists a constant $T$ such that
\begin{equation}
\left\{\begin{array}{ll}\label{3.X}
\gamma(s)=\gamma(s+T),\\
\dot{\gamma}(s)=\dot{\gamma}(s+T).
\end{array}\right.
\end{equation}
By a direct calculation, for $\gamma(s)$ in (\ref{3.111}), (\ref{3.X}) holds if and only if $\frac{b}{1-a}$ and $\frac{c}{1-a}$ are both rational numbers. Namely, $\gamma(s)$ in (\ref{3.111}) is closed if and only if $\frac{b}{1-a}$ and $\frac{c}{1-a}$ are both rational numbers.
\begin{rema}
When $b=c=0$, the geodesic in (\ref{3.111}) reduces to the expression in \cite{RO}. It is still a great circle and its length is $\frac{2\pi}{1-a}$.
\end{rema}
In other cases, the closed geodesics may not be great circles (See Fig.1. When $a=c=0$, $b=\frac{1}{2}$, the length $l=4\pi$).
\begin{figure}[H]
  \centering
  \includegraphics[width=0.4\textwidth]{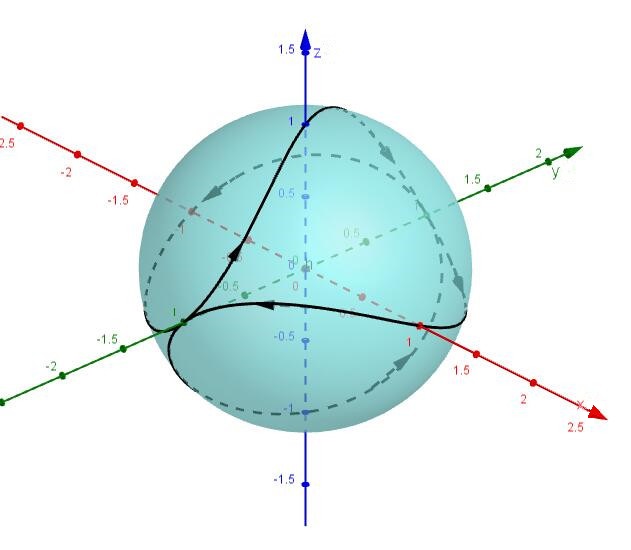}
  \caption{$a=c=0$, $b=\frac{1}{2}$, $s\in[0,4\pi]$}
\end{figure}
\begin{rema}
At point $(0,1,0)$ and $(0,-1,0)$, the geodesic in Fig.1 is self-intersecting and self-tangent, but the tangent direction is opposite. Namely, in Finsler geometry, there exist two tangent geodesics through a point, which are different from Riemannian geometry.
\end{rema}
\end{case}
\begin{case}
If there exists an irrational number among $\frac{b}{1-a}$ or $\frac{c}{1-a}$, the geodesics on $(\mathbb{S}^2,F_{Q})$ cover the sphere irregularly (See Fig.2. When $a=c=0$, $b=1-\frac{1}{\sqrt{2}}$, $s\in[0,29\pi]$), they are all non-closed.
\begin{figure}[H]
  \centering
  \includegraphics[width=0.4\textwidth]{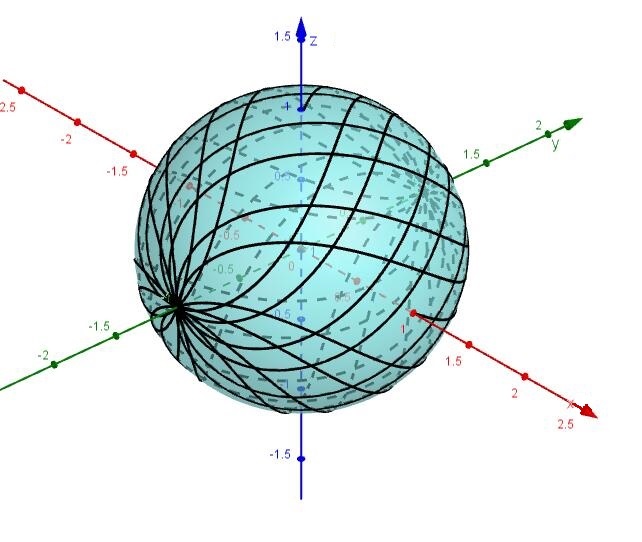}
  \caption{$a=c=0$, $b=1-\frac{1}{\sqrt{2}}$, $s\in[0,29\pi]$}
\end{figure}
\end{case}
\end{exam}
\section{Isoparametric Functions}
\subsection{Normal geodesics and focal submanifolds }
~~~~Let $\tau:M \rightarrow N$ be an embedded submanifold. For simplicity, we will denote $d\tau X$ by $X$. Define
$$\mathcal{V}(M)=\{(x,\xi)~|~x\in M,\xi\in T_x^{*}N,\xi (X)=0,\forall X\in T_xM\},$$
which is called the
{\it normal bundle} of $\tau$ or $M$. Let
$$ \mathcal{N}M={\mathcal L}^{-1}(\mathcal{V}(M))=\{(x, \eta)~|~x\in \tau(M), \eta={\mathcal L}^{-1}(\xi), \xi \in \mathcal{V}_{x}(M)\}.$$
Then $\mathcal{N}M\subset TN$. Moreover, denote
$$\mathcal{N}^{0}M=\{(x, \textbf{n})~|~x\in M, \textbf{n}\in\mathcal{N}M, F(\textbf{n})=1\}.$$
We call $\textbf{n}\in\mathcal{N}^{0}_{x}M$ the \textit{unit normal vector} of $M$.

Let $\textmd{Exp}: TN\rightarrow N$ be the exponential map of $N$. For $\eta\in\mathcal{N}_{x}M$, define the normal exponential map $E(x,\eta)=\textmd{Exp}_{x}\eta$. The focal points of $M$ are the critical values of the normal exponential map $E$. Let $M$ be an oriented hypersurface and $\textbf{n}\in\mathcal{N}^{0}_{x}M$. Set
$$\tau_{s}(x):=\textmd{Exp}(x,s\textbf{n})=\textmd{Exp}_{x}s\textbf{n}.$$
By Theorem \ref{thm0} and the notion of navigation, we immediately get the explicit expression of normal geodesics of $M$ on $(\mathbb{S}^n,F_{Q})$.
\begin{lemm}\label{coro35}
The unit speed normal geodesic $\gamma_{x}(s)$ on $(\mathbb{S}^n,F_{Q})$ with initial point $x\in M$ and initial tangent vector $\mathbf{n}$ can be expressed by
\begin{equation}\label{3.3.5.3}
\gamma(s)=\exp (sQ)((\cos s)x+\sin s(\mathbf{n}(x)-V(x))),~~~~s\in \mathbb{R}.
\end{equation}
\end{lemm}
From Lemma \ref{coro35}, we know that
\begin{equation}\label{3.3.5.3.3}
\tau_{s}(x)=\exp (sQ)((\cos s)x+\sin s(\mathbf{n}(x)-V(x))),~~ x\in M,
\end{equation}for any $s\in \mathbb{R}$.
\begin{lemm} \label{lemm23}~\cite{HCY}
Let $\tau:M  \rightarrow N(c) $ be an immersion submanifold, and let $\mathbf{n}$ be a unit normal vector to $\phi(M )$ at ${x}$.
Then $p=E(x, s\mathbf{n})$ is a focal point
 of $(M ,x)$ of multiplicity $m_0>0$ if and only if there is an eigenvalue $\lambda$ of the shape operator $A_{\mathbf{n}}$ of multiplicity $m_0$ such that
\begin{align} \lambda=\left\{
\begin{array}{rcl}
\frac{1}{s},     &      & {c=0,}\\
\sqrt{c}\cot \sqrt{c}s,   &      & {c>0,}\\
\sqrt{-c}\coth \sqrt{-c}s,       &      & {c<0.}
\end{array} \right.\label{2.3}\end{align}
\end{lemm}

Let $M$ be a connected, oriented isoparametric hypersurface in $(\mathbb{S}^n,F_{Q})$ with $\textbf{n}\in\mathcal{N}^{0}M$ and $g$ be distinct constant principal curvatures, say
\begin{align*}
\lambda_{1},~\lambda_{2},~\cdots, ~\lambda_{g}.
\end{align*}
We denote the multiplicity and the corresponding principal foliation of $\lambda_{i}$ by $m_{i}$ and $V_{i}$, respectively. Let $s_{i}=\arccot\lambda_{i}$ and $X\in V_{i}$. Set $\overline{\tau}_{s}(x)=(\cos s)x+\sin s(\textbf{n}(x)-V(x))$, then $\tau_{s}(x)=\exp(sQ)\overline{\tau}_{s}(x)$. Differentiating $\tau_{s}(x)$ in the direction $X$, from (3.41) in \cite{TP}, we have
$$\tau_{s*}X=\exp(sQ)\overline{\tau}_{s*}X=\exp(sQ)(\cos sI-\sin s\overline{A})X=\exp(sQ)\frac{\sin(s_{i}-s)}{\sin s_{i}}X,$$
where $\overline{A}$ is the shape operator of $M$ in $(\mathbb{S}^n,h)$ and on the right side we are identifying $X$ with its Euclidean parallel translate at $\tau_{s*}(x)$. $\tau_{s*}$ is injective on $V_{i}$, unless $s=s_{i}$ (mod $\pi$) for some $i$, that is, $\tau_{s_{i}*}V_{i}=0$. In this case, $\tau_{s_{i}}(x)$ is a focal point of $(M, x)$. Combine Lemma \ref{coro35}, we have
\begin{lemm}\label{coro34}
The parallel hypersurfaces and focal submanifolds of the isoparametric hypersurface $M$ in $(\mathbb{S}^n,F_{Q})$ can be expressed as $M_{s}:=\tau_{s}M,~s\neq s_{i}$ and  $M_{s_{i}}:=\tau_{s_{i}}M$, respectively, where $\tau_{s}$ is expressed by (\ref{3.3.5.3.3}).
\end{lemm}
\subsection{Proof of Theorem \ref{thm01}}
~~~~Let $\overline{f}$ be a global isoparametric function on $(\mathbb{S}^n,h)$. Set $|\nabla^{h}\overline{f}|=a(\overline{f})$, $\overline{f}(\mathbb{S}^n)=[c,d]$. From \cite{QM}, we know that $a(t)=0$ if and only if $t=c$ or $d$. In this case, $\overline{f}^{-1}(c)$ and $\overline{f}^{-1}(d)$ are only two focal submanifolds. From \cite{QM}, we know $\int^{d}_{c}\frac{1}{a(t)}dt$ is convergence. Define
$$\zeta(t)=\int_{t_{0}}^{t}\frac{1}{a(\theta)}d\theta$$
for some given $t_{0}\in(c,d)$ and any $t\in [c,d]$. Set
$$\overline{\rho}(x)=\int_{t_{0}}^{\overline{f}(x)}\frac{1}{a(t)}dt=\zeta(\overline{f}(x)),\ \ \ \ x\in\mathbb{S}^n.$$
$\overline{\rho}$ is a normalized isoparametric function with respect to $h$. Since $\zeta'(t)=\frac{1}{a(t)}>0$ when $t\in(c,d)$, $\zeta(t)$ is strictly monotonous increasing on $[c,d]$. Set $\zeta[c,d]=[\alpha,\beta]$, then
$$\overline{M}_{s}=\overline{\rho}^{-1}(s),~~ \forall s\in[\alpha,\beta],~~~~M=\overline{M}_{0}=\overline{f}^{-1}(t_0).$$

Define the map $\psi: \mathbb{S}^n\rightarrow \mathbb{S}^n$ such that \begin{equation}\label{4.5.1.1}\psi\mid_{\overline{M}_{s}}=\psi_{s}=\exp (sQ),~~~~s\in[\alpha,\beta].\end{equation}
From \cite{XM1}, we know that $\psi$ is an orientation preserving local diffeomorphism around $M$.
If we define the function $\rho$ such that $$\rho^{-1}(s)=\psi(\overline{M}_{s})=\psi(\overline{\rho}^{-1}(s)),$$
then by Lemma \ref{lemm46}, $\rho$ is a local normalized isoparametric function around $M$ on $(\mathbb S^n, F_{Q})$.

Next, we try to extend the local isoparametric function to the global one. Let $M$ be an oriented hypersurface with a given unit normal vector field $\textbf{n}$. Since $M$ is complete, we can define a distance function $\rho(x)$. If $x$ is on the side of $M$ pointed by $\textbf{n}$, then $\rho(x)=d(M,x)=\inf\{d(z,x)\mid z\in M\}$, where $d(M,x)$ is a distance from $M$ to $x$. If $x$ is on the other side of $M$, then   $\rho(x)=-d(x,M)=-\inf\{d(x,z)\mid z\in M\}$, where $d(x,M)$ is a distance from $x$ to $M$.

\begin{lemm}\label{lemm0.4}
There exists a minimal normal geodesic $\gamma(s)$ of $M$ such that $\gamma(\rho(x))=x$.
\end{lemm}
\proof Due to the definition of the distance function $\rho(x)=d(M,x)$, there exist $\{z_{n}\}\in M$ such that $s_{0}=\rho(x)=d(M ,x)=\lim\limits_{n\to+\infty}d(z_{n},x)$. Since $\mathbb{S}^n$ is compact, there exists a subsequence of $\{z_{n}\}$ which converges at $z_{0}\in\mathbb{S}^n$. Since $M$ is complete, $z_{0}\in M$ and $d(z_{0},x)=\lim\limits_{n\to+\infty}d(z_{n},x)=\rho(x)=s_{0}$.

Since $\mathbb{S}^n$ is complete and connected, there exists a minimal geodesic $\gamma:[0,s_{0}]\rightarrow \mathbb{S}^n$ satisfying $\gamma(0)=z_{0}$, $\gamma(s_{0})=x$. Then from the first variation formula of the arc-length, we know that $\dot{\gamma}(0)=\textbf{n}(z_{0})$. Namely, $\gamma(s)=E(z_{0},s\textbf{n}(z_{0}))$ is a minimal normal geodesic of $M$.
If $\rho(x)=-d(x,M)$, we can prove similarly.
\endproof

Lemma \ref{lemm0.4} also holds for $(\mathbb{S}^{n},h)$. Similarly, we can define a distance function $\overline{\rho}(x)=\overline{d}(M, x)$ on $(\mathbb S^n, h)$.

 \begin{lemm}\label{lemm0.5}
$\psi(x)=\exp(Q\int^{\overline{f}(x)}_{t_{0}}\frac{1}{a(t)}dt)x$ is a homeomorphism from $\mathbb{S}^n$ to $\mathbb{S}^n$.
\end{lemm}
\proof If there exist two different points $x_{1}$, $x_{2}\in \mathbb{S}^n$ such that $\psi x_{1}=\psi x_{2}$, then there exist two  minimal normal geodesics $\gamma_{1}(s)$ and $\gamma_{2}(s)$ of $M$ on $(\mathbb S^n, h)$ such that $x_{1}=\gamma_{1}(s_{1})$, $x_{2}=\gamma_{2}(s_{2})$. Due to the definition of $\psi$, $\psi_{s_{1}}x_{1}=\psi x_{1}=\psi x_{2}=\psi_{s_{2}}x_{2}$. Obviously,  $s_{1}\neq s_{2}$. Assume $s_{1}<s_{2}$. Then $x_{1}=\psi_{s_{2}-s_{1}}x_{2}$, which means $x_{1}$ and $x_{2}$ are on the same orbit by the one-parameter group.
Set $$\sigma(t)=\psi_{t}x_{1}=\exp (tQ)x_{1},~~~~t\in [0,s_{2}-s_{1}],$$ where $\sigma(0)=x_{1}$, $\sigma(s_{2}-s_{1})=x_{2}$. Since $|Qx|<1$, the length  $L_{\sigma}<s_{2}-s_{1}$. Hence,
 $$\overline{d}(x_{1},x_{2})\leq L_{\sigma}<s_{2}-s_{1}.$$
 Set $\overline{\gamma}_{1}(0)=z_{1}$, $\overline{\gamma}_{2}(0)=z_{2}$. Due to $$s_{2}>s_{1} +\overline{d}(x_{2},x_{1})=\overline{d}(z_{1},x_{1})+\overline{d}(x_{2},x_{1})\geq \overline{d}(z_{1},x_{2})\geq s_{2},$$
  there exits a contradiction. Hence, $\psi$ is an injection.

For any $x\in\mathbb{S}^n$, set $s_0=\rho(x)$. From Lemma \ref{lemm0.4}, there exists a minimal normal geodesic $\gamma(s)$ of $M$ on $(\mathbb S^n, F_{Q})$ such that $\gamma(s_0)=x$. Set $\gamma(0)=z_{0}\in M$ and $\dot{\gamma}(0)=\textbf{n}(z_{0})$. There also exists a minimal normal geodesic $\overline{\gamma}(s)=\exp (-sQ){\gamma}(s)$ of $M$ on $(\mathbb S^n, h)$ such that $\overline{\gamma}(0)=z_{0}$, $\dot{\overline{\gamma}}(0)=\textbf{n}(z_{0})-V$. Set $\overline{x}=\overline{\gamma}(s_0)\in \overline{M}_{s_0}$, then
$$x=\gamma(s_{0})=\exp (s_{0}Q)\overline{\gamma}(s_{0})=\exp (s_{0}Q)(\overline{x})=\psi(\overline{x}).$$
 Namely, $\psi$ is a surjection. Hence, $\psi: \mathbb{S}^{n}\rightarrow\mathbb{S}^{n}$ is a one-to-one correspondence. In addition, due to the expressions of $\psi$, $\psi$ is continuous. Since $\psi$ is a one-to-one continuous mapping, $\psi^{-1}$ is existing and continuous. This completes the proof of Lemma \ref{lemm0.5}.
\endproof
From (\ref{4.5.1.1}), we have
\begin{equation}\label{4.5.1}
\psi(x)=\exp(\overline{\rho}(x)Q)x,\ \ \ \ \ \ x\in \mathbb S^n.
\end{equation}
Combine Lemma \ref{lemm0.4} and Lemma \ref{lemm0.5}, we immediately know
\begin{equation}\label{4.4.1}
\rho(x)=\overline{\rho}(\psi^{-1}x),\ \ \ \ \ \  x\in \mathbb S^n.
\end{equation}
Set $\hat{U}=\mathbb{S}^{n}_{\overline{f}}$, $U=\psi\hat{U}$. We know $\overline{\rho}$ and $\psi$ are both smooth on $\hat{U}$. Set $P=\exp(\overline{\rho}(x)Q)\in O(n+1)$. By (\ref{4.5.1}),
$$D\psi=P(I+Qxd\overline{\rho}).$$
Due to $|d\overline{\rho}|=1$ on $\hat{U}$ and there exists a constant $q<1$ such that $|Qx|\leq q$ on $\mathbb S^n$, the Jacobian matrix $D\psi$ is invertible. Then we have
 \begin{lemm}\label{lemm0.6}
$\psi:\hat{U}\rightarrow U$ is a $C^{\infty}$ diffeomorphism.
\end{lemm}
Hence, $\rho$ is a normalized isoparametric function on $(U,F_{Q})$. Define
$f(x)=\zeta^{-1}\rho(x).$
 That is
\begin{equation}\label{4.6.2}
f(x)=\zeta^{-1}\overline{\rho}(\psi^{-1}x)=\overline{f}(\psi^{-1}x),\ \ \ \ \ \ \ x\in \mathbb{S}^n,
\end{equation}
where
\begin{equation}\label{4.6.2.1}
\psi(x)=\exp(\zeta( \overline{f}(x))Q)x=\exp(\left(\int^{\overline{f}(x)}_{t_{0}}\frac{1}{a(t)}dt\right)Q)x,\ \ \ \ \ \ \ x\in \mathbb{S}^n.
\end{equation}
Then  $f$ is smooth on $U$.

Since $f$ and $\overline{f}$ are global isoparametric functions on $(U,F_{Q})$ and $(\hat{U},h)$, respectively, they can be viewed as homogeneous functions of degree 0 on $\mathbb{R}^{n+1}$. Then $\psi$ is a homogeneous map of degree 1 on $\mathbb{R}^{n+1}$. By (\ref{4.6.2}), we can get $d{f}=d\overline fD\psi^{-1}$. That is,
\begin{equation}\label{4.7.1}
d\overline{f}=dfD\psi=dfP\left(I+Qxd\overline{\rho}\right).
\end{equation}
Since $D\psi$ and $D\psi^{-1}$ are both bounded in $U$ and $\hat{U}$, $d\overline{f}\rightarrow0$ if and only if $df\rightarrow0$. That means $f$ is $C^1$ at $x$ satisfying $df(x)=0$. Hence, by Lemma \ref{lemm0.5} and Lemma \ref{lemm0.6}, combine (\ref{4.6.2}) and (\ref{4.7.1}), we know that $f$ is a global defined isoparametric function on $(\mathbb S^n, F_{Q})$.

Conversely, if $f$ is an isoparametric function on $(\mathbb{S}^n, F_{Q})$ satisfying $F_Q(\nabla{f})=a({f})$, where $a^{2}({f})\in C^{2}(f(\mathbb S^n))$ and ${f}(\mathbb{S}^n)=[c,d]$. Then for $d(M,x)=\inf\{d(z,x)\mid z\in M\}$ (or $d(x,M)=-\inf\{d(x,z)\mid z\in M\}$) on $(\mathbb S^n, F_{Q})$, we have
\begin{lemm}\label{lemm1.11}
If $d$ is the only critical value of $f$ in $[c,d]$, then
\begin{equation}\label{1.11}
d(M,M_{d})=\inf\lim\limits_{c_{i}\rightarrow d^{-}\atop i\rightarrow \infty}d(M,M_{c_{i}})=\int_{t_{0}}^{d}\frac{1}{a(f)}df
\end{equation}
and the improper integral in (\ref{1.11}) converges.
\end{lemm}
\proof Let $\{y\}$ be any point in $M_{d}$. Since $\{f\geqslant d\}$ in $\mathbb S^n$ is closed and $\mathbb S^n$ is connected, then there are points $\{p_{i}\}\in \{f<d\}$ in $\mathbb S^n$ which tends to $y$. Namely, $\lim\limits_{i\rightarrow\infty}p_{i}=y$. Set $f(p_{i})=c_{i}$.

Let $\sigma$ be any piecewise $C^{1}$ curve which goes from $M$ to $y$. For every $c_{i}\in(c,d)$, $\sigma_{i}$ is a curve from $M$ to $M_{c_{i}}$ and $\widetilde{\sigma}_{i}$ is a curve from $M_{c_{i}}$ to $y$. Then we have
\begin{equation}\label{1.12}
\inf\limits_{\sigma}L(\sigma)\geq\inf\limits_{\widetilde{\sigma}_{i}}L(\widetilde{\sigma}_{i})+\inf\limits_{{\sigma}_{i}}L(\sigma_{i})=\inf\limits_{\widetilde{\sigma}_{i}}L(\widetilde{\sigma}_{i})+d(M,M_{c_{i}}).
\end{equation}
Since $\lim\limits_{i\rightarrow \infty}d(p_{i},y)=0$, we know $\lim\limits_{i\rightarrow \infty}\inf\limits_{\widetilde{\sigma}_{i}}L(\widetilde{\sigma}_{i})=0$. From (\ref{1.12}), we have
\begin{equation}\label{1.13}
d(M,M_{d})\geq\lim\limits_{c_{i}\rightarrow d^{-}\atop i\rightarrow \infty}d(M,M_{c_{i}}).
\end{equation}

On the other hand, since $M$ is complete, from Lemma \ref{lemm0.4} and \cite{HYS}, we can find $\{x_{i}\}\in M$ and normal geodesics $\gamma_{i}$ such that $L(\gamma_{i})=d(x_{i},p_{i})=d(M, M_{c_{i}})$. Due to $d(x_{i},p_{i})+d(p_{i},y)\geq d(x_{i},y)\geq d(M,M_{d})$, we have
\begin{equation}\label{1.14}
\lim\limits_{i\rightarrow \infty}L(\gamma_{i})=\lim\limits_{c_{i}\rightarrow d^{-}\atop i\rightarrow \infty}d(M,M_{c_{i}})\geq d(M,M_{d}).
\end{equation}
From \cite{HYS}, combine (\ref{1.13}) and (\ref{1.14}),
$$d(M,M_{d})=\lim\limits_{c_{i}\rightarrow d^{-}\atop i\rightarrow \infty}d(M,M_{c_{i}})=\int_{t_{0}}^{d}\frac{1}{a(f)}df.$$
\endproof
From Lemma \ref{lemm1.11}, $\int^{d}_{c}\frac{1}{a(t)}dt$ is convergence. In addition, for $a^{2}(t)\in C^{2}(f(\mathbb{S}^{n}))$, using similar method in \cite{QM}, we have
\begin{lemm}\label{lemm1.12}
$f$ has no critical value in $(c,d)$.
\end{lemm}
Then from Lemma \ref{lemm1.11} and Lemma \ref{lemm1.12}, similarly, we can define $\widetilde{\psi}(x)=\exp(-\zeta( {f}(x))Q)x$ and prove $\overline{f}=f\circ\widetilde{\psi}^{-1}$ is an isoparametric function on $(\mathbb{S}^n, h)$. This completes the proof of Theorem \ref{thm01}.
\subsection{Homogeneous functions and homogeneous polynomials}
~~~~For Riemannian isoparametric functions, we have known the following fact from \cite{TP}. Let $\phi: \mathbb{R}^{n+1}\rightarrow \mathbb{R}$ be a homogeneous function. Then $f=\phi\mid_{(\mathbb S^{n},h)}$ is an isoparametric function on $(\mathbb S^n,h)$ if and only if there exist a smooth function $\widehat{a}(\phi,r)$ and a continuous function $\widehat{b}(\phi,r)$ such that $\phi$ satisfies
\begin{equation}\label{1.2}
\left\{\begin{array}{ll}
|\nabla^E\phi|^2=\widehat{a}(\phi,r),\\
\Delta^{E}\phi=\widehat{b}(\phi,r),
\end{array}\right.
\end{equation}
where $r=|x|$.

We try to deduce similar equations to characterize isoparametric functions on a Randers sphere. Recall that a Randers metric $F$ can be defined by navigation from the datum $(h, V)$, where $h=\sqrt{h_{ij}y^{i}y^{j}}$, $V=v^{i}\frac{\partial}{\partial x^{i}}$. Set $V_{0}=v_{i}y^{i}=h_{ij}v^{j}y^{i}$, $\lambda=1-\|V\|_{h}^{2}$.
\begin{lemm} \label{lemm41}~\cite{HYS1}
Let $(N, F, d\mu_{BH})$ be an $n$-dimensional Randers space, where $F=\frac{\sqrt{\lambda h^{2}+V_{0}^{2}}-V_{0}}{\lambda}$. Then $f$ is an isoparametric function if and only if $f$ satisfies
$$\left\{\begin{array}{ll}
|df|_h+\langle df, V^*\rangle_h=\tilde a(f),\\
\frac{1}{ |df|_h}\Delta^hf+\text{div}V +\frac{1}{ |df|_h^2}\langle d\langle df,V^*\rangle_h,df\rangle_h=\frac{\tilde b(f)}{\tilde a(f)}.
\end{array}\right.$$
\end{lemm}
When $V$ is a Killing vector field, we have
$$\text{div}_{\sigma} V=h^{ij}v_{i|j}=0.$$
\begin{lemm} \label{lemm42}~\cite{TP}
Let~$(\mathbb S^n,h)\hookrightarrow \mathbb R^{n+1} (n\geq2)$ be the standard Euclidean sphere and $\phi: \mathbb R^{n+1}\rightarrow \mathbb R$ be a homogeneous function of degree $k$. Then
\begin{align}\label{na}
\nabla^h\phi=\nabla^E\phi-k\phi x,
\end{align}
\begin{align}\label{na1}
\Delta^h\phi=\Delta^E\phi-k(k+n-1)\phi,
\end{align}
where $\nabla^E$ and $\Delta^E$ denote the Euclidean gradient and  Laplacian in $\mathbb R^{n+1}$, respectively.
\end{lemm}
If $\phi$ is a positively homogeneous function of degree $k$, (\ref{na}) and (\ref{na1}) still holds. By Lemma \ref{lemm41} and Lemma \ref{lemm42}, we get the following
\begin{theo}\label{thm41}
Let $\phi: \mathbb R^{n+1}\rightarrow \mathbb R$ be a positively homogeneous function of degree $k$. Then $f=\phi\mid_{\mathbb S^{n}}$ is an isoparametric function on $(\mathbb S^n, F_{Q})$ if and only if $\phi$ satisfies
\begin{equation}\label{1.3}
\left\{\begin{array}{ll}
|\nabla^E\phi-k\phi x|+\langle \nabla^E\phi,xQ\rangle=\tilde a(\phi,r),\\
\frac{\Delta^E\phi-k(k+n-1)\phi}{ |\nabla^E\phi-k\phi x|}
      +\frac{\langle \nabla^E\langle \nabla^E\phi,xQ\rangle,\nabla^E\phi\rangle}{|\nabla^E\phi-k\phi x|^2}=\tilde b(\phi,r).
\end{array}\right.
\end{equation}
\end{theo}
Theorem \ref{thm41} gives equations to characterize isoparametric functions on a Randers sphere. But the equations are hard to solve. Hence, we use Theorem \ref{thm01} to derive the solution of (\ref{1.3}).

Let $\overline{\phi}(x)$ be a homogeneous function of degree $k$ on $\mathbb R^{n+1}$ which satisfies (\ref{1.2}). Then $\overline{f}=\overline{\phi}\mid_{\mathbb{S}^n}$ is a global defined isoparametric function on $(\mathbb{S}^n,h)$ and $\overline{\phi}(x)=|x|^{k}\overline{f}(\frac{x}{|x|})$. By Theorem \ref{thm01}, $f=\overline{f}\circ\psi^{-1}$ is a global defined isoparametric function on $(\mathbb{S}^n, F_{Q})$. Set $\phi=|x|^{k}f(\frac{x}{|x|})$. $\phi$ is a positively homogeneous function of degree $k$ on $\mathbb R^{n+1}$ such that $f=\phi\mid_{\mathbb{S}^n}$, which must be a solution of (\ref{1.3}).

Set $x=\psi(\overline{x})=\exp (\zeta (\overline{f}(\overline{x}))Q)\overline{x}$. Obviously, $|x|=|\overline{x}|$. Then
\begin{equation}\label{4.8}
\phi(x)=|x|^{k}f(\frac{x}{|x|})=|\overline x|^{k}\overline{f}\circ\psi^{-1}(\frac{x}{|\overline{x}|}),
\ \ \ \ \ \ x\in\mathbb{R}^{n+1}.\end{equation}
Since $\psi^{-1}$ is a homogeneous map of degree 1,
\begin{equation}\label{4.8.2}
\phi(x)=|\overline x|^{k}\overline f(\frac{\overline{x}}{|\overline{x}|})=\overline{\phi}(\overline{x})
=\overline{\phi}(\psi^{-1}x),
\ \ \ \ \ \ x\in\mathbb{R}^{n+1}.\end{equation}
From Lemma \ref{lemm42}, $|\nabla^h\overline{\phi}|^2=|\nabla^E\overline{\phi}|^2-k^2\overline{\phi}^2,$ that is, $a(t)^2=\widehat{a}
(t,1)-k^{2}t^{2}$.  Thus, we have
\begin{equation}\label{4.8.1}
\psi(\overline{x})=\exp\left(Q\int^{\frac{\overline{\phi}(\overline x)}{|\overline{x}|^{k}}}_{t_{0}}\frac{1}{\sqrt{\widehat{a}
(t,1)-k^{2}t^{2}}}dt\right)\overline x, \ \ \ \ \ \ \overline x\in\mathbb{R}^{n+1}.
\end{equation}
Combine (\ref{4.8.2}) and (\ref{4.8.1}), we give the solution of (\ref{1.3}).

Conversly, if there exists a positively homogeneous function $\phi$ of degree $k$ on $\mathbb R^{n+1}$ which satisfies (\ref{1.3}), then $f=\phi|_{\mathbb{S}^n}$ is an isoparametric function on $(\mathbb{S}^n, F_{Q})$ and $\phi=|x|^{k}f(\frac{x}{|x|})$. From Theorem \ref{thm01}, $\overline{f}=f\circ\widetilde{\psi}^{-1}$ is an isoparametric function on $(\mathbb{S}^n, h)$, where $\widetilde{\psi}(x)=\exp(-\zeta( {f}(x))Q)x$. We can also find a homogeneous function $\overline{\phi}$ of degree $k$ which satisfies $\overline{\phi}\mid_{\mathbb{S}^{n}}=\overline{f}$ and $\overline{\phi}=\phi\circ\widetilde{\psi}^{-1}$. Namely, we have the following
\begin{theo}$\overline{\phi}$ is a homogeneous function of degree $k$ satisfying (\ref{1.2}) if and only if $\phi=\overline{\phi}\circ\psi^{-1}$ is a positively homogeneous function of degree $k$ satisfying (\ref{1.3}), where $\psi$ is expressed by (\ref{4.8.1}).
\end{theo}
\begin{lemm} \label{lemm021}~\cite{HDY}
Let~$M$ be an anisotropic submanifold in a Randers space~$(N,F,d\mu_{BH})$ with the navigation datum~$(h,V)$. If~$F$ has isotropic~$\mathbf{S}$-curvature~$\mathbf{S}=(n+1)c(x)F$, then for any $\mathbf{n} \in \mathcal{N}^0(M)$, the shape operators of $M$ in Randers space~$(N, F)$ and Riemannian space~$(N, h)$, ${A}_{\mathbf{n}}$ and $\bar{A}_{\bar{\mathbf{n}}}$,  have the same principal vectors and satisfy
\begin{align}\label{3.24}\lambda=\bar{\lambda}+c(x),\end{align}
where~$\lambda$ and~$\bar{\lambda}$ are the principal curvatures of~$M$ in a Randers space~$(N, F)$ and a Riemannian space~$(N, h)$, respectively.\end{lemm}
From Lemma \ref{lemm021}, the number and multiplicity of principal curvatures are the same by navigation. Hence, there exists the same isoparametric hypersurface in $(\mathbb{S}^n,h)$ and $(\mathbb{S}^n,F_{Q})$. So from Theorem A, we immediately have
\begin{prop}
Let $M$ be a connected isoparametric hypersurface in $(\mathbb S^{n},F_{Q})$ with $g$ principal curvatures $\lambda_{i}=cot\theta_{i}$, $0<\theta_{i}<\pi$, with multiplicities $m_{i}$. Then $M$ is an open subset of a level set of the restriction to $\mathbb S^{n}$ of a Cartan-M\"{u}nzner polynomial $\phi$ of degree $g$.
\end{prop}
\subsection{Proof of Theorem \ref{thm02}}
~~~~For a given connected isoparametric family $M_{t}$ in $(\mathbb{S}^n,F_{Q})$, there exists an isoparametric function $\tilde{f}:(\mathbb{S}^n,F_{Q})\rightarrow[c,d]$ such that $M_{t}\subset \widetilde{f}^{-1}(t)$ and $M=M_{t_{0}}\subseteq \tilde{f}^{-1}(t_{0})=\tilde{f}^{-1}(\frac{1}{2}(c+d))$. Let $\rho$ be a normalized isoparametric function such that $\nabla \rho=\frac{\nabla \tilde{f}}{F_Q(\nabla \tilde{f})}$ and $\rho(\mathbb{S}^n_{\tilde{f}})=(-\alpha,\alpha)$, where $\alpha>0$ is a constant. Set $\rho^{-1}(0)=\tilde{f}^{-1}(t_{0})$, then $M\subset\rho^{-1}(0)$. From  Lemma \ref{lemm021}, $M$ is also a connected isoparametric hypersurface in $(\mathbb{S}^n,h)$. From Theorem A, there exists a Cartan-M\"{u}nzner polynomial $\phi$ of degree $g$ such that $\overline{f}=\phi\mid_{\mathbb{S}^n}$ is an isoparametric function on $(\mathbb{S}^n,h)$ and $M$ is an open subset of a level set of $\overline{f}$.

In this case, from (\ref{1.4}) and Lemma \ref{lemm42}, $$a(\overline{f})=|\nabla^{h}\phi|=\sqrt{|\nabla^{E}\phi|^{2}-g^2\phi^2}\big|_{\mathbb{S}^{n}}=\sqrt{g^{2}-g^2\overline{f}^2}=g\sqrt{1-\overline{f}^2}.$$
Hence, $a(t)=g\sqrt{1-t^2}$. Obviously, $\zeta(t)=\int_{{0}}^{t}\frac{1}{a(\theta)}d\theta=\frac{1}{g}\arcsin t$. From (\ref{4.6.2.1}), we have
$$\psi(x)=\exp(\frac{1}{g}\arcsin \overline{f}(x)Q)x=\exp(\frac{1}{g}\arcsin\frac{\phi(x)}{|x|^{g}}Q)x.$$
Due to Theorem \ref{thm01}, $f=\overline{f}\circ\psi^{-1}$ is an isoparametric function on $(\mathbb{S}^n,F_{Q})$ and each $M_{t}$ is an open subset of a level set of $f$.
Namely, there exists a $t'$ such that $M_{t'}\subset f^{-1}(t)$ and $M=f^{-1}(0)$.  Meanwhile, from Lemma \ref{lemm32} and Lemma \ref{lemm0.4},
\begin{equation}\label{0.01}
x=\exp(sQ)\overline{x}=\exp(\overline{d}(M,\overline{x})Q)\overline{x}.
\end{equation}
Hence, from (\ref{3.3.5.3.3}),
\begin{align}\label{0.001}
&f^{-1}(t)=E(x,\zeta(t)\textbf{n})=\tau_{\zeta(t)}(x)=\exp(\zeta(t)Q)\overline{\tau}_{\zeta(t)}(x)\nonumber\\
&=\{\exp(\frac{\arcsin t}{g}Q)((\cos\frac{\arcsin t}{g})x+\sin(\frac{\arcsin t}{g})(\textbf{n}-V))\mid x\in M\},-1<t<1.
\end{align}

Conversely, for a given Cartan-M\"{u}nzner polynomial $\phi$ of degree $g$, set $f=\phi\circ\psi^{-1}\mid_{\mathbb{S}^n}$. Then $f^{-1}(t)$ is a connected isoparametric family in  $(\mathbb{S}^n,F_{Q})$, where $-1<t<1$. Namely, we can find an isoparametric family in $(\mathbb{S}^n,F_{Q})$ through a positively homogeneous polynomial, which is equivalent to construct a mean curvature flow.

Further more, we discuss the focal submanifolds in $(\mathbb S^n, F_{Q})$. By Theorem B, there exist only two connected focal submaifolds $\overline{M}_{+}=\overline{\tau}_{\frac{\pi}{2g}}(x)$ and $\overline{M}_{-}=\overline{\tau}_{-\frac{\pi}{2g}}(x)$ in $(\mathbb S^n, h)$. Let $\overline{x}$ be a focal point of $M$ in $(\mathbb S^n, h)$. From (\ref{4.7.1}), $\nabla \overline{f}(\overline{x})\rightarrow0$ if and only if $\nabla f(x)\rightarrow0$. Hence, $x$ is a focal point of $M$ in $(\mathbb S^n, F_{Q})$. By (\ref{0.01}), there also exist two focal submanifolds $M_{\pm}$ in $(\mathbb S^n, F_{Q})$ corresponding to $\overline{M}_{\pm}$ in $(\mathbb S^n, h)$. From (\ref{0.001}), we have
\begin{equation}\label{4.7.2}
M_{\pm}=\tau_{\pm \frac{\pi}{2g}}(x)=\{\exp(\pm \frac{\pi}{2g}Q)(\cos\frac{\pi}{2g}x\pm \sin\frac{\pi}{2g}(\textbf{n}-V))\mid x\in M\}.
\end{equation}
Theorem \ref{thm02} can be proved.
\subsection{Examples of isoparametric functions and focal submanifolds}
~~~~In this section, we will give some examples of isoparametric functions on $(\mathbb{S}^n,F_{Q})$, the corresponding isoparametric families and focal submanifolds. We are familiar with the fact that the standard form of real antisymmetric matrices by the orthogonal similarity transformation can be written as
$$Q=\left(
\begin{array}{cccccccccc}
Q_{1}& & & & & & \\
     &Q_{2}& & & & & \\
     & &\ddots& & & & \\
     & & &Q_{j}& & & \\
     & & & &0& & \\
     & & & & &\ddots& \\
     & & & & & &0\\
\end{array}
\right),$$
where \begin{equation}
Q_{i}=\left(
\begin{array}{cc}
0&a_{i}\\
-a_{i}&0
\end{array}
\right), \ \ \ \ \ \ |a_{i}|<1, \ \ \ \ \ \ 1\leq i\leq j.
\end{equation}
If $\overline{\phi}$ is a homogeneous polynomial of degree $g$, then $\phi$ is a positively homogeneous function of degree $g$.  Using (\ref{4.8.2}), we have $\phi(x)=\overline{\phi}(\psi^{-1}x)$, where
\begin{equation}\label{4.12}
\psi(x)=(\sum\limits_{i=1}^{j}I_{i}\cos a_{i}\left(\frac{1}{g}\arcsin\frac{\phi(x)}{|x|^{g}}\right)+\frac{1}{a_{i}}P_{i}\sin a_{i}\left(\frac{1}{g}\arcsin\frac{\phi(x)}{|x|^{g}}\right))x+I_{n+1-2j}x,
\end{equation}
where $I_{i}$ is a $(n+1)\times(n+1)$ matrix which satisfies each row and column of $I_{i}$ is 0, except for the elements of row $2i-1$, column $2i-1$ and row $2i$, column $2i$ are both $1$. $P_{i}$ is a $(n+1)\times(n+1)$ matrix which satisfies each row and column of $P_{i}$ is 0, except for the element of row $2i$, column $2i-1$ is $-1$ and row $2i-1$, column $2i$ is $1$. And $I_{n+1-2j}$ is still a $(n+1)\times(n+1)$ matrix which satisfies the last $n+1-2j$ elements on the diagonal of $I_{n+1-2j}$ are all 1, all the other elements are 0.
\begin{exam} \label{exam0}
$\textbf{(g=1)}$ Let $(\mathbb{S}^n, F_{Q})$ be an $n$-dimensional Randers sphere, where $F=\frac{\sqrt{\lambda h^{2}+V_{0}^{2}}-V_{0}}{\lambda}$. Then $\overline{f}$ is an isoparametric function if and only if there exists a homogeneous polynomial $\overline{\phi}:\mathbb{R}^{n+1}\rightarrow\mathbb{R}$ which satisfies
$$\overline{\phi}(x)=\langle x,p\rangle,\ \ \ \ p\in\mathbb{R}^{n+1},\ \ \ \ |p|=1,$$
such that $\overline{f}=\overline{\phi}\mid_{\mathbb{S}^n}$.
Then we have $\phi(x)=\overline{\phi}(\psi^{-1}x)$. Using (\ref{4.12}),
\begin{equation*}
\psi(x)=(\sum\limits_{i=1}^{j}I_{i}\cos a_{i}(\arcsin\frac{\langle x,p\rangle}{|x|})+\frac{1}{a_{i}}P_{i}\sin a_{i}(\arcsin\frac{\langle x,p\rangle}{|x|}))x+I_{n+1-2j}x.
\end{equation*}
In Riemannian case, for this linear height function, consider the regular level sets $\overline{M}_{s}=\{x\in\mathbb{S}^{n}\mid\langle x,p\rangle=\sin s\}$, $-\frac{\pi}{2}<s<\frac{\pi}{2}$. The two focal submanifolds are $\overline{M}_{+}=M_{\frac{\pi}{2}}=\{p\}$ and $\overline{M}_{-}=M_{-\frac{\pi}{2}}=\{-p\}$, respectively. From (\ref{0.001}) and (\ref{4.7.2}), the regular level sets $$M_{s}=\{x\in\mathbb{S}^{n}\mid\langle x,\exp(sQ)p\rangle=\sin s\},\ \ \ \  -\frac{\pi}{2}<s<\frac{\pi}{2}.$$
The connected two focal submanifolds $M_{+}=\exp(\frac{\pi}{2}Q)\{p\}$, $M_{-}=\exp(-\frac{\pi}{2}Q)\{-p\}$ are still two points.
\end{exam}
\begin{rema} \label{rema03}
Take
$Q=\left(
\begin{array}{cc}
Q'&0\\
0&0\\
\end{array}
\right)$
and $p=e_{n+1}$, where $Q'\in o(n)$. Example \ref{exam0} reduces to Example $\textit{3.11}$ in~\cite{HDY}.
\end{rema}
Take $p=(1,0,0)$ and
\begin{equation*}
Q=\left(
\begin{array}{cccc}
0&0&\frac{1}{2}\\
0&0&0\\
-\frac{1}{2}&0&0
\end{array}
\right),
\end{equation*}
we can give the image of mean curvature flow in $(\mathbb{S}^{2},h)$ and $(\mathbb{S}^{2},F_{Q})$ (See Fig. 3 and Fig. 4). In this case, the two focal submanifolds
$M_{+}=\{(\frac{\sqrt{2}}{2},0,-\frac{\sqrt{2}}{2})\}$ and $M_{-}=\{(-\frac{\sqrt{2}}{2},0,-\frac{\sqrt{2}}{2})\}$ are two points.
\begin{figure}[H]
\centering
    \begin{tabular}{cc}
        \begin{minipage}[t]{3.5in}
        \includegraphics[width=2.5in]{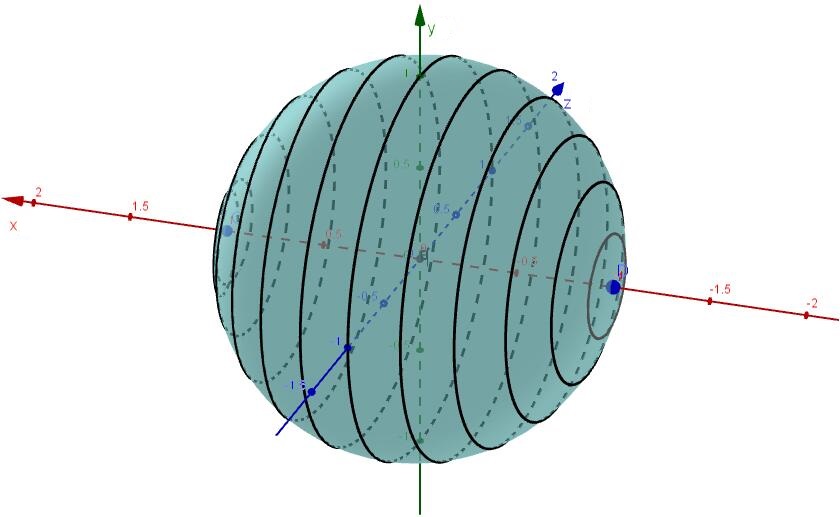}
        \caption{Isoparametric family in $(\mathbb{S}^{2},h)$}
        \end{minipage}
        \begin{minipage}[t]{3.5in}
        \includegraphics[width=2.5in]{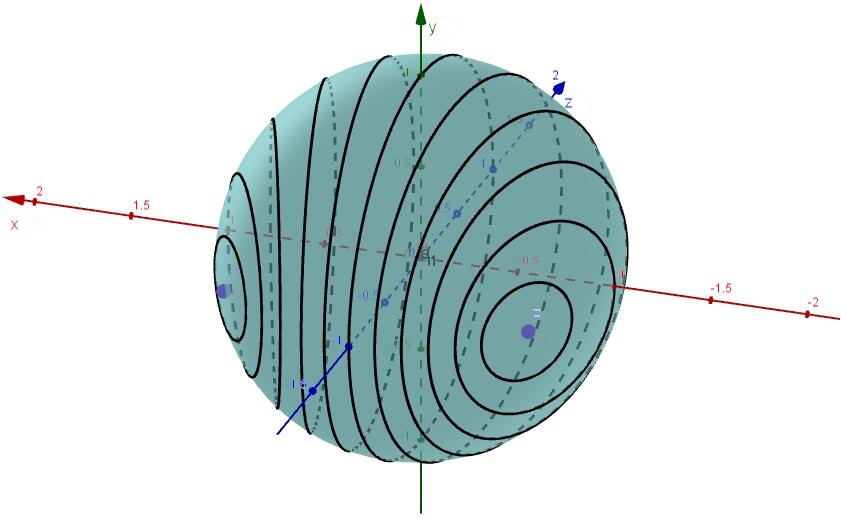}
        \caption{Isoparametric family in $(\mathbb{S}^{2},F_{Q})$}
        \end{minipage}
    \end{tabular}
\end{figure}
\begin{exam} \label{exam01}
$\textbf{(g=2)}$ Define the function $\overline{\phi}:\mathbb{R}^{p+1}\times \mathbb{R}^{q+1}\rightarrow \mathbb{R}$
$$\ \ \ \ \ \ \ \  \ \ \ \ \ \ \ \ \ \ \ \ \ \ \ \ \ \ \ \ \ \ \ \ \ \ \ \ \ \ \ \ (\overline{x}_{1}\ \ ,\ \ \overline{x}_{2})\ \ \ \ \mapsto |\overline{x}_{1}|^2-|\overline{x}_{2}|^{2},$$
where $\overline{x}_{1}=(x_{1}, x_{2},\dots, x_{p+1})$, $\overline{x}_{2}=(x_{p+2}, x_{p+3}\dots, x_{p+q+2})$, $p,q\in N^{+}$ and $p+q=n-1$. $\overline{\phi}$ is a homogeneous polynomial.

Set $Q=(q^{\alpha}_{\beta})_{(n+1)\times(n+1)}$, where $q^{p+2}_{p+1}=-q^{p+1}_{p+2}=a$ and the others are $0$. By (\ref{4.12}), we obtain $\phi(x)=\phi(\psi^{-1}x)$, where
\begin{align*}
\psi(x)&=(x_{1},\dots,x_{p},\cos (a\arcsin\frac{|\overline{x}_{1}|^{2}-|\overline{x}_{2}|^{2}}{|\overline{x}_{1}|^{2}+|\overline{x}_{2}|^{2}})x_{p+1}+sin (a\arcsin\frac{|\overline{x}_{1}|^{2}-|\overline{x}_{2}|^{2}}{|\overline{x}_{1}|^{2}+|\overline{x}_{2}|^{2}})x_{p+2},\\
&-\sin(a\arcsin\frac{|\overline{x}_{1}|^{2}-|\overline{x}_{2}|^{2}}{|\overline{x}_{1}|^{2}+|\overline{x}_{2}|^{2}})x_{p+1}+\cos (a\arcsin\frac{|\overline{x}_{1}|^{2}-|\overline{x}_{2}|^{2}}{|\overline{x}_{1}|^{2}+|\overline{x}_{2}|^{2}})x_{p+2},x_{p+3},\dots,x_{n+1}).
\end{align*}
In Riemannian case, consider the regular level sets $\overline{M}_{s}=\{x=(\overline{x}_{1},\overline{x}_{2})\in\mathbb{S}^{n}\mid|\overline{x}_{1}|^2-|\overline{x}_{2}|^{2}=\sin 2s\}$, $-\frac{\pi}{4}<s<\frac{\pi}{4}$. The two focal submanifolds are $\overline{M}_{+}=\overline{M}_{\frac{\pi}{4}}=S^{p}\times\{0\}$ and $\overline{M}_{-}=\overline{M}_{-\frac{\pi}{4}}=\{0\}\times\mathbb{S}^{q}$, respectively. From (\ref{0.001}) and (\ref{4.7.2}), the regular level sets can be expressed as
\begin{align*}
M_{s}=\{(&x_{1},\dots,x_{p},\cos(as)x_{p+1}+\sin(as)x_{p+2},-\sin(as)x_{p+1}+\cos(as)x_{p+2},\\
&x_{p+3},\dots,x_{n+1})\in\mathbb{S}^n\mid|\overline{x}_{1}|^2-|\overline{x}_{2}|^{2}=\sin 2s\},\ \ \ \  -\frac{\pi}{4}<s<\frac{\pi}{4}.
\end{align*}
The two focal submanifolds are
\begin{align*}
M_{+}&=(\exp \frac{\pi}{4}Q)\{S^{p}\times\{0\}\}\\
     &=\{(x_{1},\dots,x_{p+2},0,\dots,0)\mid x_{p+2}+\tan(\frac{\pi}{4}a)x_{p+1}=0,\sum\limits_{i=1}^{p+2}x_{i}^{2}=1\}
\end{align*}
and
\begin{align*}
M_{-}&=\exp(-\frac{\pi}{4}Q)\{\{0\}\times\mathbb{S}^{q}\}\\
     &=\{(0,\dots,0,x_{p+1},\dots,x_{n+1})\mid x_{p+1}+\tan(\frac{\pi}{4}a)x_{p+2}=0,\sum\limits_{i=p+1}^{n+1}x_{i}^{2}=1\}.
\end{align*}
Namely, $M_{+}$ is $\mathbb{S}^{p}$ in $\mathbb{R}^{p+2}$ and $M_{-}$ is $\mathbb{S}^{q}$ in $\mathbb{R}^{q+2}$.
\end{exam}

Yali Chen \\
School of Mathematical Sciences, Tongji University, Shanghai, 200092, China\\
E-mail: chenyl90@tongji.edu.cn\\

Qun He \\
School of Mathematical Sciences, Tongji University, Shanghai, 200092, China\\
E-mail: hequn@tongji.edu.cn\\

\end{document}